\documentclass[12pt,reqno]{amsart}

\usepackage{graphics}
\usepackage{graphicx}
\usepackage{color}

\usepackage{amsfonts}
\usepackage{eufrak}

\usepackage{amsmath}
\usepackage{amstext}
\usepackage{amsopn}
\usepackage{amsbsy}
\usepackage{amscd}
\usepackage{amsxtra}
\usepackage{amsthm}

\usepackage{fancyhdr}
\newcommand{\tz}{\tilde{z}}
\newcommand{\ta}{\tilde{a}}

\newcommand{\tV}{\tilde{V}}

\newcommand{\tl}{\tilde{\lambda}}

\newcommand{\ds}{\displaystyle}

\def \R {{\bf R}}

\def \a {\alpha}
\def \b {\beta}
\def \D {\Delta}

\def \l {\lambda}

\def \O {\Omega}
\def \Om {\Omega}

\def \p {\partial}
\def \s {\sigma}

\def \ve {\varepsilon}
\def \i {\infty}
\def \ds {\displaystyle}

\def \bd {\bigtriangledown}

\def \h {\hspace{.5cm}}

\def \h* {\hspace*{1cm}}

\newcommand{\beq}{\begin{equation}}
\newcommand{\eeq}{\end{equation}}
\newcommand{\pp}{\partial}

\newtheorem{lemma}{Lemma}[section]

\newtheorem{theorem}[lemma]{Theorem}

\newtheorem{remark}{Remark}
\newtheorem{claim}{Claim}

\begin{document}
\title[Symbiotic solitons in coupled NLS ]
{ Symbiotic bright solitary wave solutions of coupled nonlinear
Schr\"{o}dinger equations}
%\title{On  Coupled Nonlinear Schr\"{o}dinger Equations in Bose-Einstein
%Condensates}

\author{Tai-Chia Lin}
\address{Department of Mathematics,
  National Taiwan University,
  Taipei 106, Taiwan,} \address{
  National Center of Theoretical Sciences, National Tsing
Hua University, Hsinchu, Taiwan}\email{tclin@math.ntu.edu.tw}

\author{Juncheng Wei}
\address{Department of Mathematics,
  The Chinese University of Hong Kong,
  Shatin, Hong Kong }
\email{wei@math.cuhk.edu.hk}


\date{}

\maketitle

\begin{abstract}

Conventionally, bright solitary wave solutions can be obtained in
self-focusing nonlinear Schr\"{o}dinger equations with attractive
self-interaction. However, when self-interaction becomes
repulsive, it seems impossible to have bright solitary wave
solution. Here we show that there exists symbiotic bright solitary
wave solution of coupled nonlinear Schr\"{o}dinger equations with
repulsive self-interaction but strongly attractive interspecies
interaction. For such coupled nonlinear Schr\"{o}dinger equations
in two and three dimensional domains, we prove the existence of
least energy solutions and study   the location and configuration
of symbiotic bright solitons. We use Nehari's manifold to
construct least energy solutions and derive their asymptotic
behaviors by some techniques of singular perturbation problems.
\end{abstract}

\section{Introduction}

In this paper, we study symbiotic bright solitary wave solutions
of two-component system of time-dependent nonlinear
Schr\"{o}dinger equations called Gross-Pitaevskii equations given
by \beq\label{1.3}\left\{\begin{array}{l}
 i \hbar\pp_t\psi_1=-\frac{\hbar^2}{2m}\Delta\psi_1+\tilde{V}_1(x)\psi_1+U_{11}|\psi_1|^2\psi_1+U_{12}|\psi_2|^2\psi_1,\\
 i \hbar\pp_t\psi_2=-\frac{\hbar^2}{2m}\Delta\psi_2+\tilde{V}_2(x)\psi_2+U_{22}|\psi_2|^2\psi_2+U_{12}|\psi_1|^2\psi_2,\ x\in\Omega,\ t>0.\ \ \ \ \
\end{array}\right.
\eeq which models a binary mixture of Bose-Einstein condensates
with two different hyperfine states called a double condensate.
Here $\O\subseteq \mathbb{R}^N (N\leq 3)$ is the domain for
condensate dwelling, $\psi_j$'s are corresponding condensate wave
functions, $\hbar$ is the Planck constant divided by $2\pi$ and
$m$ is atom mass. The constants $U_{jj}\sim a_{jj}$, $j=1$,$2$,
and $U_{12}\sim a_{12}$, where $a_{jj}$ is the intraspecies
scattering length of the $j$-th hyperfine state and $a_{12}$ is
the interspecies scattering length. Besides, $\tilde{V}_j$ is the
trapping potential for the $j$-th hyperfine state. In physics, the
usual trapping potential is given by $$ \tV_j(x)= \sum_{k=1}^N
\ta_{j,k}(x_k-\tz_{j,k})^2\quad\hbox{ for
}\:x=(x_1,\cdots,x_N)\in\Omega, j=1,2\,,$$ where $\ta_{j,k}\geq 0$
is the associated axial frequency, and
$\tz_j=(\tz_{j,1},\cdots,\tz_{j,N})$ is the center of the trapping
potential $\tV_j$.

When the constant $U_{jj}$ is negative and large enough,
self-interaction of the $j$-th hyperfine state is strongly
attractive and the associated condensate tends to increase its
density at the centre of the trap potential in order to lower the
interaction energy (cf.~\cite{21}). This may result in spikes and
bright solitons which can be observed experimentally in three
dimensional domain (cf.~\cite{CTW}). Conversely, when the constant
$U_{jj}$ becomes positive, self-interaction on the $j$-th
hyperfine state turns into repulsion which cannot support the
existence of bright solitons. To create bright solitons while each
self-repulsive state cannot support a soliton by itself, the
interspecies attraction may open a way to make two-component
solitons called symbiotic bright solitons. Recently, symbiotic
bright solitons in only one dimensional domain have been
investigated as the interspecies scattering length $a_{12}$ is
negative and sufficiently large (cf.~\cite{PB}). However, in two
and three dimensional domains, the existence of symbiotic bright
solitons has not yet been proved. In this paper, we want to show
the existence of such solitons by studying the least energy
solutions of two-component system of nonlinear Schr\"{o}dinger
equations.

To obtain symbiotic bright solitons in a double condensate, we may
set $\psi_1(x,t)= u(x)\,e^{i\,\tl_1\,t}$, $\psi_2(x,t)=
v(x)\,e^{i\,\tl_2\,t}$ and use Feshbach resonance to let
$U_{jj}$'s, $\tl_j$'s and $\ta_{j,k}$'s be very large quantities.
By rescaling and some simple assumptions, the system (\ref{1.3})
with very large $U_{jj}$'s, $\tl_j$'s and $\ta_{j,k}$'s is
equivalent to the following singularly perturbed problem:
\begin{equation}{\label{eq:1-1}}
\left\{\begin{array}{l}
\ve^2 \D u- V_1 (x) u+\mu_1 u^3 + \b u v^2 =0\ \ \ \mbox{in}\ \ \O,\\
\ve^2 \D v- V_2 (x)  v+\mu_2 v^3 + \b u^2 v =0\ \ \ \mbox{in}\ \ \O,\\
u,v >0\ \ \mbox{in}\ \ \O,\\
u=v=0\ \ \mbox{on}\ \ \p\O,
\end{array}\right.
\end{equation}
where $u$ and $v$ are corresponding condensate amplitudes, $\ve>0$
is a small parameter, and $\b\sim -a_{12}\neq 0$ is a coupling
constant. Here we may use the zero Dirichlet boundary condition
which may come from~\cite{GP}. To study symbiotic bright solitons
of double condensates, we consider two cases of the domain $\O$.
One is to set $\O$ as the entire space $\R^N (N\leq 3)$. The other
is to set $\O$ as a bounded smooth domain in $\R^N$. The constants
$\mu_j\sim -U_{jj}\leq 0\,, j=1, 2\,,$ give repulsive
self-interaction, and $\b\sim -a_{12}>0$ means attractive
interaction of solutions $u$ and $v$. Moreover, $V_j>0\,, j=1,2$
are the associated trapping potentials.

Another motivation of studying the problem (\ref{eq:1-1}) may come
from the formation of bright solitons in a mixture of a degenerate
Fermi gas with a Bose-Einstein condensate in the presence of a
sufficiently attractive boson-fermion interaction. Recently, there
have been successful observations and associated experimental and
theoretical studies of mixtures of a degenerate Fermi gas and a
Bose-Einstein condensate (cf.~\cite{DJ}, \cite{MRR} and \cite{M}).
Recently, the corresponding model has been given by
\beq\label{bf1} \left\{\begin{array}{l}
 i \hbar\pp_t\varphi^B=-\frac{\hbar^2}{2m_B}\Delta\varphi^B+ V_B(x)\varphi^B
 +\,g_B N_B|\varphi^B|^2\varphi^B+ g_{BF}\,\ds\sum_{j=1}^{N_F}\,|\varphi^F_j|^2\varphi^B\,,\\
 i
 \hbar\pp_t\varphi^F_j=-\frac{\hbar^2}{2m_F}\Delta\varphi^F_j+ V_F(x)\varphi^F_j
 + g_{BF}\,N_B\,|\varphi^B|^2\varphi^F_j\,,\: x\in\Omega,\ t>0\,,
 j=1,\cdots, N_F\,,
\end{array}\right.
\eeq where $N_B$ and $N_F$ are the numbers, $m_B$ and $m_F$ are
the mass of bosons and fermions, $V_B$ and $V_F$ are trap
potentials, $\varphi^B$ and $\varphi^F_j$'s are wave functions of
Bose-Einstein condensate and individual fermions, respectively.
When the constant $g_B$ is positive i.e. repulsive
self-interaction, and the constant $g_{BF}$ is negative and large
enough enough i.e. strongly attractive interspecies interaction,
bright solitons may appear in such a system. Using the
system~(\ref{bf1}) (cf.~\cite{KB}), a novel scheme to realize
bright solitons in one-dimensional atomic quantum gases (i.e. the
domain $\Omega$ is one dimensional) can be found. Here we want to
study bright solitons in two and three-dimensional atomic quantum
gases i.e. the domain $\Omega$ is of two and three dimensional. As
for the problem~(\ref{eq:1-1}), we may set $\varphi^B=
u(x)\,e^{i\,\tl_1\,t}/\sqrt{N_B}$, $\varphi^F_j=
v_j(x)\,e^{i\,\tl_2\,t}$ and suitable scales on $m_B, m_F, V_B,
V_F, g_B, g_{BF}$ and $\tl_j$'s. Then the system~(\ref{bf1}) can
be transformed into
\begin{equation}{\label{bf2}} \left\{\begin{array}{l}
\ve^2 \D u- V_1 (x) u+\mu_1 u^3 + \b u \ds\sum_{j=1}^{N_F}\,v_j^2 =0\ \ \ \mbox{in}\ \ \O,\\
\ve^2 \D v_j- V_2 (x)  v_j + \b u^2 v_j =0\ \ \ \mbox{in}\ \ \O,\quad j=1,\cdots,N_F\,,\\
u,v_j >0\ \ \mbox{in}\ \ \O,\\
u=v=0\ \ \mbox{on}\ \ \p\O\,,
\end{array}\right.
\end{equation} which can be generalized as a singular perturbation
problem given by
\begin{equation}{\label{bf3}} \left\{\begin{array}{l}
\ve^2 \D u- V_1 (x) u+\mu_1 u^3 + \b u \ds\sum_{j=1}^{m}\,v_j^2 =0\ \ \ \mbox{in}\ \ \O,\\
\ve^2 \D v_j- V_2 (x)  v_j +\mu_2 v_j^3 +  \b u^2 v_j =0\ \ \ \mbox{in}\ \ \O,\quad j=1,\cdots,m\,,\\
u,v_j >0\ \ \mbox{in}\ \ \O,\\
u=v=0\ \ \mbox{on}\ \ \p\O\,,
\end{array}\right.
\end{equation} where $\mu_j\leq 0, j=1, 2$ are constants and
$m=N_F\in \mathbb{N}$. In particular, the problem~(\ref{bf3})
becomes the problem~(\ref{eq:1-1}) as $m=1$.

In this paper, we study the asymptotic behavior of so-called
least-energy solutions of the problem~(\ref{eq:1-1}) which may
give symbiotic bright solitons in two and three dimensional
domains. By this, we mean
\begin{enumerate}
\item $(u_\ve,v_\ve)$ is a solution of (\ref{eq:1-1}), \item
$E_{\ve, \Omega, V_1, V_2} [u_\ve,v_\ve]\leq E_{\ve, \Omega, V_1,
V_2} [u,v]$ for any nontrivial solution $(u,v)$ of (\ref{eq:1-1}),
\end{enumerate}
where $E_{\ve, \Om, V_1, V_2} [u,v]$ is the energy functional
defined as follows:
\begin{eqnarray}{\label{eq:1-2}}
E_{\ve, \Om, V_1, V_2} [u,v] & :=\  & \frac{\ve^2}{2}\int_\O |\bd
u|^2 + \frac{V_1}{2} \int_\O u^2 -\frac{\mu_1}{4} \int_\O u^4\\
\nonumber & & +~~~\frac{\ve^2}{2} \int_\O |\bd v|^2 +
\frac{V_2}{2} \int_\O v^2 - \frac{\mu_2}{4} \int_\O v^4 \\
\nonumber & & -~~~\frac{\b}{2}\int_\O u^2 v^2\,,
\end{eqnarray} for $u,v\in H_0^1(\O)$. Actually, it is easy to generalize our results
to the problem~(\ref{bf3}) for $m\in\mathbb{N}$. In the case of
$\O=\R^N, N=2,3$, the least energy solution is also called ground
state. In our previous papers~\cite{lw1}, \cite{lw2} and
\cite{lw3},  we studied the existence and asymptotics of least
energy solutions when $ \mu_1$ and $\mu_2$ are positive constants.
Hereafter, we study the case that both $\mu_1$ and $\mu_2$ are
non-positive constants.

As $\b \leq \sqrt{\mu_1 \mu_2}$, it is obvious that
\begin{equation}
\label{ncon} \int_\O [ \ve^2 |\nabla u|^2 + V_1 u^2 + \ve^2
|\nabla v|^2 + V_2 v^2] =  \int_{\O} [2 \b u^2 v^2  +\mu_1  u^4 +
\mu_2 v^4] \leq 0
\end{equation}
for any $(u, v)$ satisfying the problem~(\ref{eq:1-1}) and hence $
u, v \equiv 0$. To get nontrivial solutions of the
problem~(\ref{eq:1-1}), the assumption $ \b> \sqrt{ \mu_1 \mu_2}$
is necessary. So throughout the paper, we assume that
\begin{equation}
\label{basic} \mu_1\leq 0,\quad \mu_2\leq 0,\quad \b > \sqrt{\mu_1
\mu_2}\,.
\end{equation}

To study least energy solutions, we define a Nehari manifold
\begin{equation}
\label{Nehari} N(\ve, \O, V_1, V_2)= \Biggl\{  (u,v)\in H_0^1 (\O)
\times H_0^1(\O)  \Biggl| \begin{array}{l} \int_\O [\ve^2 |\bd
u|^2 +V_1  u^2 +
\ve^2  |\bd v|^2 + V_2  v^2]\\
 =  \int_\O [ 2 \b u^2 v^2 +\mu_1 u^4 +\mu_2 v^4]
\end{array}
 \Biggl\}.
\end{equation}
Note that here, unlike \cite{lw1}-\cite{lw3}, the Nehari manifold
$N(\ve, \O, V_1, V_2)$ has only one constraint. On such a
manifold, we consider the minimization problem given by
\begin{equation}{\label{eq:2-2}}
c_{\ve, \O, V_1, V_2} :=\  \inf_{(u,v)\in  N(\ve,\O, V_1, V_2),
\atop{u, v\geq 0, \atop{u, v\not\equiv 0}}} E_{\ve, \Omega, V_1,
V_2} [u,v]\,.
\end{equation}
When $\ve=1$, $V_j\equiv \lambda_j>0, j=1,2$ i.e. constant
trapping potentials and the domain $\O=\R^N$, the Euler-Lagrange
equations of the problem~(\ref{eq:2-2}) are
\begin{equation}{\label{eq:1-4g}}
\left\{ \begin{array}{l}
\D u -\l_1 u +\mu_1 u^3 + \b u v^2 =\  0 \ \ \mbox{in}\ \ \R^N,\\
\D v -\l_2 v +\mu_2 v^3 + \b u^2 v =\  0 \ \ \mbox{in}\ \ \R^N,\\
u, v  \to 0 \ \ \mbox{as}\ \ |y|\to +\i.
\end{array}\right.
\end{equation}
For such a problem, we have
\begin{theorem}
\label{t1.1} Assume that (\ref{basic}) holds. Then $ c_{1, \R^N,
\l_1, \l_2}$ is attained and hence the problem~(\ref{eq:1-4g})
admits a ground state solution which is radially symmetric and
strictly decreasing.
\end{theorem}

Now we consider the existence of ground state solutions for
nonconstant trapping potentials. Namely, we consider the problem
of coupled nonlinear Schr\"{o}dinger equations given by
\begin{equation}{\label{eq:1-4s}}
\left\{ \begin{array}{l}
\ve^2 \D u -V_1 (x) u +\mu_1 u^3 + \b u v^2 =\  0 \ \ \mbox{in}\ \R^N,\\
\ve^2 \D v -V_2 (x) v +\mu_2 v^3 + \b u^2 v =\  0 \ \ \mbox{in}\ \R^N,\\
u, v  \to 0 \ \ \mbox{as}\ \ |y|\to +\i,
\end{array}\right.
\end{equation}
where $V_j$'s satisfy
\begin{equation}
 0<b_j^0= \inf_{x \in \R^N} V_j(x) \leq  \lim_{|x| \to \infty} V_j(x)= b_j^\infty \leq +\infty,
 \quad j=1,2\,.
\end{equation}
Then we have the following theorem on the existence of ground
state solutions of the problem~(\ref{eq:1-4s}).
\begin{theorem}
\label{t1.2}
If either $b_1^\infty+ b_2^\infty = +\infty$ or
\begin{equation}
\label{c7} c_{\ve, \R^N, V_1, V_2} <  c_{\ve, \R^N, b_1^\infty,
b_2^\infty}
\end{equation}
Then $ c_{\ve, \R^N, V_1, V_2}$ is attained and hence the
problem~(\ref{eq:1-4s}) admits a ground state solution.
\end{theorem}
\noindent Our next theorem is to show the asymptotic behavior of
these ground state solutions as follows:
\begin{theorem}\label{t1.3}
Assume (\ref{basic}) and
\begin{equation}\label{rho1}
\inf_{x\in\mathbb{R}^n} c_{1, \R^N, V_1(x),V_2(x) }< c_{1, \R^N,  b_1^\infty, b_2^\infty}
\end{equation} hold. Then
\begin{itemize}
\item[(i)]$ c_{\ve, \R^N, V_1, V_2}$ is attained and the
problem~(\ref{eq:1-4s}) admits a ground state solution $(u_\ve,
v_\ve)$.

\item[(ii)]Let $P^\ve$  and $Q^\ve$ be the unique local maximum
points of $u_{\ve}$ and $ v_\ve$ respectively. Let
$u_{\ve}(P^\ve+\ve y):=U_{\ve}(y), v_{\ve} (Q^\ve +\ve y):=V_\ve
(y) $. Then as $\ve\rightarrow0$, $(U_{\ve},V_{\ve})\rightarrow(U,
V)$, where $(U, V)$ satisfies (\ref{eq:1-4g}). Furthermore,
\begin{equation}\label{p2}
\frac{|P^\ve- Q^\ve|}{\ve}\rightarrow0\,,\quad c_{1, \R^N,
V_1(P^\ve), V_2 (Q^\ve)}  \to \inf_{x \in \R^N} c_{1, \R^N, V_1
(x), V_2 (x)}\,.
\end{equation}
\end{itemize}
\end{theorem}

\begin{remark}~~In general, the condition~(\ref{rho1}) is difficult to check.
However, if $\inf\limits_{x\in\R^N}\,V_j
(x)<\lim\limits_{|x|\rightarrow+\infty}\,V_j(x), j=1,2$, then
(\ref{rho1}) is satisfied. \end{remark}

Theorem~\ref{t1.3} can be extended to general bounded domains.
Firstly, we set $\O$ as a bounded smooth domain and trapping
potentials $V_j$'s as constants $\lambda_j$'s. Namely, we consider
the following system
\begin{equation}{\label{eq:1-1d}}
\left\{\begin{array}{l}
\ve^2 \D u- \l_1  u+\mu_1 u^3 + \b u v^2 =0\ \ \ \mbox{in}\ \ \O,\\
\ve^2 \D v- \l_2  v+\mu_2 v^3 + \b u^2 v =0\ \ \ \mbox{in}\ \ \O,\\
u,v >0\ \ \mbox{in}\ \ \O,\\
u=v=0\ \ \mbox{on}\ \ \p\O\,.
\end{array}\right.
\end{equation}
The asymptotic behavior of corresponding least energy solutions
can be characterized by
\begin{theorem}\label{t1.4}~~For any $\b >\sqrt{\mu_1 \mu_2}$ and $\ve$
sufficiently small, the problem~(\ref{eq:1-1d}) has a least energy
solution $(u_\ve,v_\ve)$. Let $P_\ve$ and $Q_\ve$ be the local
maximum points of $u_\ve$ and $v_\ve$, respectively. Then
 $|P_\ve-Q_\ve|/\ve \to 0$,
\begin{eqnarray}{\label{eq:1-3}}
d(P_\ve,\p\O)\to \max_{P\in \O} d(P,\p\O),\  d(Q_\ve,\p\O) \to
\max_{P\in \O} d(P,\p\O)\,,
\end{eqnarray}
and $u_\ve(x), v_\ve(x)\to 0$ in $C_{loc}^1 (\bar\O\backslash
 \{P_\ve,
Q_\ve\})$. Furthermore, as $\ve \to 0$, $(U_\ve,V_\ve)\to
(U_0,V_0)$ which is a least-energy solution of~(\ref{eq:1-4g}),
where
$$U_\ve(y):=\ u_\ve(P_\ve+\ve y),\quad V_\ve(y):=\ v_\ve(P_\ve+\ve y)\,.$$

\end{theorem}
\noindent By Theorem~\ref{t1.4}, we may generalize
Theorem~\ref{t1.3} to bounded smooth domains. The main idea may
follow the proof of Corollary~2.7 in~\cite{lw3}. Moreover, by the
same arguments of Theorems~\ref{t1.1}-\ref{t1.4}, one may get
similar results for the problem~(\ref{bf3}).

As $\mu_1, \mu_2>0$, the assumption $\beta<\beta_0$ is essential
in our previous works (cf.~\cite{lw1}-\cite{lw3}) for the
existence and the asymptotic behaviors of ground state (least
energy) solutions, where $0<\beta_0<\sqrt{\mu_1\,\mu_2}$ is a
small constant. For larger $\beta$'s, results of ground and bound
state solutions can be found in \cite{ac1}, \cite{BWW}, \cite{P1}
and \cite{S1}. On the other hand, when the sign of $\mu_j$'s
becomes negative i.e.~$\mu_1, \mu_2\leq 0$, the assumption of
$\beta$'s can be changed as $\beta>\sqrt{\mu_1\,\mu_2}$ which is
sufficient to prove the existence and the asymptotic behaviors of
ground state solutions (see~Theorem~\ref{t1.1}-\ref{t1.4}). These
are new results of two and three dimensional bright solitary wave
solutions for negative $\mu_j$'s.

Conventionally, there has been a vast literature on the study of
concentration phenomena for single singularly perturbed nonlinear
Schr\"{o}dinger equations with attractive self-interaction. See
\cite{amn}, \cite{BDS}, \cite{bw1}, \cite{cl}, \cite{df},
\cite{df1}, \cite{df2}, \cite{DY}, \cite{GW2}, \cite{jt},
\cite{kw}, \cite{Li}, \cite{mm}, \cite{W1}, \cite{WW1},
\cite{zwang} and the references therein. In particular, a good
survey can be found in \cite{n1} and \cite{n2}. However, until
now, there are only few papers working on systems of coupled
nonlinear Schr\"{o}dinger equations, especially for two and three
dimensional Bose-Einstein condensates. This paper seems to be the
first in showing rigorously that strong interspecies attraction
may produce symbiotic bright solitons in two and three dimensional
Bose-Einstein condensates even though self-interactions are
repulsive.

The organization of this paper is as follows:

\noindent In Section~2, we extend the classical Nehari's manifold
approach to a system of semilinear elliptic equations in order to
find a least energy solution to~ the problem~(\ref{eq:1-1}).
Hereafter, we need the condition $\beta>\sqrt{\mu_1\,\mu_2}$ for
strong interspecies attraction. Using approximation argument and
energy upper bound, we may show Theorem~\ref{t1.1}, \ref{t1.2} and
Theorem \ref{t1.3} in Section~3 and 4, respectively. In Section~5,
we follow the same ideas of~\cite{lw1} to complete the proof of
Theorem~\ref{t1.4}.

Throughout this paper, unless otherwise stated, the letter $C$
will always denote various generic constants which are independent
of $\ve$, for $\ve$ sufficiently small. The constant $\sigma \in
(0, \frac{1}{100})$ is a fixed small constant.

\vskip 0.5cm

\noindent {\bf Acknowledgments:} The research of the first author
is partially supported by a research Grant from NSC of Taiwan. The
research of the second author is partially supported by an
Earmarked Grant from RGC of Hong Kong. The authors also want to
express their sincere thanks to the referee's suggestions.

\section{Nehari's Manifold Approach : Existence of a Least-Energy Solution
to (\ref{eq:1-1})}

\setcounter{equation}{0}

In this section, we use Nehari's manifold approach to obtain a
least energy solution to (\ref{eq:1-1}). Nehari's manifold
approach has been used successfully in the study of single
equations.  Conti et al~\cite{CTV} have used Nehari's manifold to
study solutions of competing species systems which are related to
an optimal partition problem in $N$-dimensional domains. In our
previous paper \cite{lw1}, we also used Nehari's manifold approach
to find least energy solutions and symbiotic bright solitons.

We consider the following minimization problem
\begin{equation}\label{eq:2-2n}
c_{\ve, \O, V_1, V_2} :=\  \inf_{(u,v)\in N (\ve,\O, V_1, V_2),
\atop{u, v\geq 0, \atop{u, v\not\equiv 0}}} E_{\ve, \Omega, V_1,
V_2} [u,v]
\end{equation}
where $N (\ve, \O, V_1, V_2)$ and $E_{\ve, \O, V_1, V_2}$ are
defined in Section~1. Note that, for $N\leq 3$, by the compactness
of Sobolev embedding $H_0^1(\O) \hookrightarrow L^4 (\O),
N(\ve,\O, V_1, V_2)$ and $c_{\ve, \O, V_1, V_2}$ are well-defined.
Now we want to show that
\begin{theorem}
\label{t2.1} Let $\O$ be a smooth and bounded domain in $\R^N,
N\leq 3$. Suppose that $\b >\sqrt{\mu_1\mu_2}.$ Then for $\ve$
sufficiently small, $c_{\ve, \O, V_1, V_2}$ can be attained by
some $(u_\ve,v_\ve)\in N(\ve,\O, V_1, V_2)$ satisfying
\begin{equation}
\label{uepvep}
 C_1 \ve^N  \leq  \int_\Omega u_\ve^4  \leq  C_2 \ve^N,
 \ C_1 \ve^N  \leq  \int_\Omega v_\ve^4  \leq C_2 \ve^N,
\end{equation}
where $C_1, C_2$ are two positive constants independent of $\ve $
and $\O$.
\end{theorem}

%\begin{proof}

We first note that if $(u,v)\in N(\ve,\O, V_1, V_2)$, then
\begin{eqnarray}{\label{eq:2-3}}
E_{\ve, \Omega, V_1, V_2} [u,v] & =\  & \frac 1 4 \Biggl( \ve^2
\int_\O |\bd u|^2 +
\int_\O V_1 u^2 + \ve^2 \int_\O |\bd v|^2 + \int_\O V_2 v^2\Biggl) \\
\nonumber & =\  & \frac 1 4 \Biggl[\mu_1 \int_\O u^4 + 2\b \int_\O
u^2 v^2 + \mu_2 \int_\O v^4 \Biggl].
\end{eqnarray}
\noindent Let $(u_n,v_n)$ be a minimizing sequence. Then by
Sobolev embedding $H_0^1 (\O)\hookrightarrow L^q (\O)$ for
$1<q<\frac{2N}{N-2}$, we see that $u_n\to u_\ve$, $v_n\to
v_\ve$(up to a subsequence) for some functions $u_\ve\geq 0$,
$v_\ve\geq 0$ in $L^4(\O)$ and hence
\begin{equation}{\label{eq:2-4}}
E_{\ve, \Omega, V_1, V_2} [u_n,v_n]\to \frac 1 4 \Biggl[
\mu_1\int_\O u_\ve^4 + 2\b \int_\O u_\ve^2 v_\ve^2 +\mu_2 \int_\O
v_\ve^4\Biggl]=\ c_{\ve, \O, V_1, V_2}.
\end{equation}
By (\ref{eq:2-4}) and the weak lower semicontinuity of the $H^1$
norm, we have
\begin{equation}{\label{eq:2-5}}
c_{\ve, \O, V_1, V_2} \geq \frac 1 4 \Biggl( \ve^2\int_\O |\bd
u_\ve|^2 + \int_\O V_1 u_\ve^2 +\ve^2 \int_\O |\bd v_\ve|^2 +
\int_\O V_2 v_\ve^2\Biggl),
\end{equation}
and
\begin{equation}{\label{eq:2-6}}
\ve^2 \int_\O |\bd u_\ve|^2 + \int_\O V_1 u_\ve^2  +\ve^2 \int_\O
|\bd v_\ve|^2 +  \int_\O V_2 v_\ve^2  \leq  \mu_1 \int_\O u_\ve^4
+2 \b \int_\O u_\ve^2 v_\ve^2+ \mu_2 \int_\O v_\ve^4\,.
\end{equation}

Next we consider for $t >0$,
\begin{equation}{\label{eq:2-8}}
\b_{(u, v)} (t) =\  E_{\ve, \Omega, V_1, V_2} [\sqrt{t} u,
\sqrt{t} v]\,.
\end{equation}
Our first claim is

\begin{claim}
If $ 2 \b \int_\O u^2 v^2 +\mu_1 \int_\O u^4 +\mu_2 \int_\O v^4
>0$, then $\b_{(u, v)}(t)$ attains a unique maximum point $t_0$,
where
\begin{equation}
\label{eq:2-25} t_0=\frac{ \int_\O [ \ve^2 |\bd u|^2 +V_1 u^2
+\ve^2 |\bd v|^2 + V_2 v^2]}{\int_\O [ 2\b u^2 v^2 + \mu_1 u^4
+\mu_2 v^4] }\,.
\end{equation}
Furthermore, $ (\sqrt{t_0}u, \sqrt{t_0} v) \in N (\ve, \O, V_1,
V_2)$.
\end{claim}
\begin{proof}~~Since
\begin{eqnarray*}
\b_{(u, v)} (t) & =\  & t \Biggl[\frac{\ve^2}{2} \int_\O |\bd u|^2 +
 \frac{1}{2}
\int_\O V_1 u^2 + \frac{\ve^2}{2} \int_\O |\bd v|^2 + \frac{1}{2} \int_\O V_2
v^2\Biggl]  \\
 & & -  t^2 \Biggl[ \frac{\mu_1}{4} \int_\O u^4 + \frac{\mu_2}{4} \int_\O v^4 +\frac 1 2 \b \int_\O  u^2 v^2
 \Biggl]\,,
\end{eqnarray*}
then the proof follows by simple calculations. We omit the details
here.
\end{proof}
\noindent By Claim~1 and proper choice of $(u,v)$, it is easy to
check that the Nehari manifold $N (\ve, \O, V_1, V_2)$ is
nonempty.  Our second claim is
\begin{claim}
The inequalities of (\ref{uepvep}) hold if $\beta >\sqrt{\mu_1
\mu_2}$.
\end{claim}

\vskip 0.5cm

\begin{proof}~~ We first prove the upper bound of $ c_{\ve, \O, V_1, V_2}$.
Since $ \beta >\sqrt{\mu_1 \mu_2}$, there exists $ \alpha\neq 0$
such that $ 2 \b \alpha^2  + \mu_1 \alpha +\mu_2 >0$. In fact, we
may set $\alpha=-\frac{ \mu_2}{\mu_1}$ if $\mu_j<0, j=1,2$. For
$\ve$ sufficiently small, we choose a test function $w$ such that
$ \mbox{support} ( w) \subset B_{\ve } (P)$ where $ P \in \Om$.
Let $ (u, v)= (\alpha w,  w)$. Then $ \int_\O [ 2 \b u^2 v^2
+\mu_1 u^4+\mu_2 v^4]>0$.  By Claim~1, there exists $ t_0>0$
independent of $\varepsilon$ such that $ (\sqrt{t_0} u, \sqrt{t_0}
v) \in N(\ve, \O, V_1, V_2)$. Hence we obtain
\begin{equation}
\label{cep1} c_{\ve, \O, V_1, V_2} \leq C \ve^N\,,
\end{equation} where $C$ is a positive constant independent of
$\varepsilon$ and $\O$. Combining (\ref{cep1}) with
(\ref{eq:2-3}), we obtain that \begin{equation}\label{2.9-1}
\int_\Om [ \ve^2 |\nabla u_\ve|^2 +V_1 u_\ve^2 + \ve^2 |\nabla
v_\ve^2| + V_2 v_\ve^2] \leq C_2 \ve^N\,.\end{equation} For
(\ref{2.9-1}), we may rescale spatial variables by $\varepsilon$
and apply the standard Gagliardo-Nirenberg-Sobolev inequality in
$\R^N$ (cf.~\cite{E1}). Consequently, \beq\label{ub1}\int_\O
u_\ve^4 \leq C_2 \ve^N,\quad \int_\Om v_\ve^4 \leq C_2
\ve^N\,,\eeq where $C_2$ is a positive constant independent of
$\varepsilon$ and $\O$.

For lower bound estimates, the definition of the manifold $N(\ve,
\O, V_1, V_2)$ may give
\[
\int_\Om  [ \ve^2 |\nabla u|^2 +V_1 u^2 + \ve^2 |\nabla v^2| + V_2
v^2] \leq  2 \b \int_\O u^2 v^2\,,\] for any $(u, v) \in N (\ve,
\O, V_1, V_2)$. On the other hand, as for (\ref{ub1}), we may
rescale spatial variables by $\varepsilon$ and apply the standard
Gagliardo-Nirenberg-Sobolev inequality in $\R^N$ (cf.~\cite{E1})
to derive
\[ \int_\Om  [ \ve^2 |\nabla u|^2 +V_1 u^2 + \ve^2 |\nabla v^2| + V_2 v^2]
\geq C \ve^{N/2} \left[(\int_\O u^4)^{1/2} +(\int_\O
v^4)^{1/2}\right] \geq C \ve^{N/2} (\int_\O u^2 v^2)^{1/2}
\] for any $(u, v) \in N (\ve,\O, V_1, V_2)$, and hence we obtain that for any $ (u, v) \in
N(\ve, \O, V_1, V_2)$, $(u, v) \not \equiv (0, 0)$,
\begin{equation}
\label{uv} \int_\O u^2 v^2 \geq C \ve^N\,,
\end{equation} where $C$ is a positive constant independent of
$\varepsilon$ and $\O$. Due to $\int_\O u^2 v^2\leq \left(\int_\O
u^4\right)^{1/2}\,\left(\int_\O v^4\right)^{1/2}$, (\ref{ub1}) and
(\ref{uv}) may yield lower bound estimates $ \int_\O u_\ve^4 \geq
C_1 \ve^N$ and $\int_\O v_\ve^4 \geq C_1 \ve^N$, where $C_1$ is a
positive constant independent of $\varepsilon$ and $\O$.

\end{proof}

Finally we claim that
\begin{lemma}
$(u_\ve, v_\ve)$ is a least-energy solution of (\ref{eq:1-1}).
\end{lemma}

\begin{proof}

By Claim~2 and (\ref{eq:2-6}), we have $ 2 \b \int_\O u_\ve^2
v_\ve^2 +\mu_1 \int_\O u_\ve^4 +\mu_2 \int_\O v_\ve^4 >0$.
Moreover, by Claim~1, there exists $t_0>0$ such that $(\sqrt{t_0}
u_\ve, \sqrt{t_0} v_\ve) \in N (\ve,\O, V_1, V_2)$ i.e.
\begin{equation}{\label{eq:2-12}}
\ve^2 \int_\O |\bd u_\ve|^2 +  \int_\O V_1 u_\ve^2 +\ve^2 \int_\O
|\bd v_\ve|^2 +  \int_\O V_2 v_\ve^2  = t_0 \left[  \mu_1 \int_\O
u_\ve^4 +2 \b  \int_\O u_\ve^2 v_\ve^2 +  \mu_2 \int_\O v_\ve^4
\right]\,.
\end{equation}
Consequently, (\ref{eq:2-6}) and (\ref{eq:2-12}) may give
\begin{equation}
\label{t01}  t_0 \leq 1\,.
\end{equation}
On the other hand,
\begin{equation}{\label{eq:2-16}}
E_{\ve, \Omega, V_1, V_2} [\sqrt{t_0}u_\ve, \sqrt{t_0}v_\ve] \geq
c_{\ve, \O, V_1, V_2} =\  \frac 1 4 \Biggl[ \mu_1 \int_\O u_\ve^4
+ 2\b \int_\O u_\ve^2 v_\ve^2 + \mu_2 \int_\O v_\ve^4 \Biggl],
\end{equation}

\begin{equation}{\label{eq:2-17}}
E_{\ve, \Omega, V_1, V_2} [\sqrt{t_0}u_\ve, \sqrt{t_0}v_\ve] =
t_0^2  \frac 1 4 \Biggl[ \mu_1 \int_\O u_\ve^4 + 2\b  \int_\O
u_\ve^2 v_\ve^2 + \mu_2 \int_\O v_\ve^4 \Biggl].
\end{equation}
Since $t_0>0$,  (\ref{eq:2-16}) and (\ref{eq:2-17}) imply that  $
t_0 \geq 1$. Thus by (\ref{t01}), we obtain $ t_0=1$ and
$(u_\ve,v_\ve)\in N(\ve,\O, V_1, V_2)$. Therefore, $(u_\ve,v_\ve)$
attains the minimum $c_{\ve, \O, V_1, V_2}$.

Now we want to claim that $(u_\ve, v_\ve)$ is a nontrivial
solution of (\ref{eq:1-1}). Since $(u_\ve, v_\ve)$ is an energy
minimizer on the Nehari manifold $N(\ve,\O, V_1, V_2)$, there
exists a Lagrange multiplier $\a$ such that
\begin{equation}{\label{eq:2-20}}
\bd E_{\ve, \Omega, V_1, V_2} [u_\ve,v_\ve] + \a \bd
G[u_\ve,v_\ve] =\ 0\,,
\end{equation}
where
\begin{equation}{\label{eq:2-21}}
G[u,v] =\ \int_\O  [ \ve^2 |\bd u|^2 +V_1  u^2  +\ve^2 |\bd v|^2
+V_2  v^2 ] - \int_\O[ \mu_1  u^4 +2\b  u^2 v^2+\mu_2 v^4]\,.
\end{equation}
Acting (\ref{eq:2-20}) with $(u_\ve, v_\ve)$, and making use of
the fact that $ (u_\ve, v_\ve) \in N(\ve, \O, V_1, V_2)$, we see
that
\[ \alpha \int_\O 2 [  \ve^2 |\bd u_\ve|^2 +V_1  u_\ve^2  +\ve^2 |\bd v_\ve|^2 +V_2  v_\ve^2
] - 8 \alpha\int_\O[ \mu_1  u_\ve^4 +2\b  u_\ve^2 v_\ve^2+\mu_2
v_\ve^4]=0\,,
\] and
\[ \alpha \int_\O[ \mu_1  u_\ve^4 +2\b  u_\ve^2 v_\ve^2+\mu_2 v_\ve^4]=0\,.\]
Since $ (u_\ve, v_\ve) \not \equiv (0, 0)$ and
\[ \int_\O[ \mu_1  u_\ve^4 +2\b  u_\ve^2 v_\ve^2+\mu_2 v_\ve^4]
= \int_\O  [ \ve^2 |\bd u_\ve|^2 +V_1  u_\ve^2  +\ve^2 |\bd
v_\ve|^2 +V_2  v_\ve^2] >0\,,\] then $ \alpha=0$. This proves that
$$\bd E_{\ve, \Omega, V_1, V_2} [u_\ve,v_\ve]=\ 0$$
and hence $(u_\ve,v_\ve)$ is a critical point of $E_{\ve, \Omega,
V_1, V_2} [u,v]$ and satisfies (\ref{eq:1-1}). By Hopf boundary
Lemma, it is easy to show that $u_\ve>0$ and $v_\ve>0$. Therefore,
we may complete the proof of this Lemma and Theorem~\ref{t2.1}.

\end{proof}

Another useful characterization of $c_{\ve, \O, V_1, V_2}$ is
given as follows:
\begin{lemma}
\label{cnew}
If $\b > \sqrt{\mu_1 \mu_2}$, then we have
\begin{eqnarray}\label{eq:2-25s}
c_{\ve, \O, V_1, V_2} &= &\  \inf_{u,v\in H_0^1(\O), \ u\not\equiv 0,
v\not\equiv 0, \atop{\int_\O [ 2\b u^2 v^2 +\mu_1 u^4 + \mu_2 v^4]>0} } \sup_{t>0}
E_{\ve, \O, V_1, V_2}[\sqrt{t}u, \sqrt{t}v] \\
&= &  \inf_{u,v\in H_0^1(\O), \ u\not\equiv 0, v\not\equiv 0,
\atop{\int_\O [ 2\b u^2 v^2 +\mu_1 u^4 + \mu_2 v^4]>0} } \frac{
\int_\O [ |\nabla u|^2 + V_1 u^2 + |\nabla v|^2 +V_2 v^2]}{
(\int_\O [ 2\b u^2 v^2 + \mu_1 u^4 +\mu_2 v^4])^{\frac{1}{2}}}.
\notag
\end{eqnarray}
\end{lemma}

\begin{proof}
The last identity in (\ref{eq:2-25s}) follows from simple
calculations. To prove (\ref{eq:2-25s}), we denote the right hand
side of (\ref{eq:2-25s}) by $m_\ve$. From Theorem~\ref{t2.1},
$c_{\ve, \O, V_1, V_2}$ is attained at $(u_\ve,v_\ve)\in N(\ve,\O,
V_1, V_2)$. Moreover, by Claim~1 in Theorem~\ref{t2.1}, $E_{\ve,
\Om, V_1, V_2} [\sqrt{t} u_\ve, \sqrt{t}v_\ve]$ attains its
maximum at $t=1$. Hence
\begin{equation}{\label{eq:2-26}}
m_\ve \leq c_{\ve, \O, V_1, V_2} =\  E_{\ve, \Om, V_1, V_2}
[u_\ve, v_\ve] =\  \sup_{t>0} E_{\ve, \Om, V_1, V_2} [\sqrt{t}
u_\ve, \sqrt{t} v_\ve].
\end{equation}
On the other hand, fix $u,v\in H_0^1(\O)$ such that $u, v\geq 0$
and $\int_\O [ 2 \b u^2 v^2 + \mu_1 u^4 + \mu_2 v^4]
>0$. Let $t_0$ be a critical point of $\b_{(u, v)} (t)$. Then
$(\sqrt{t_0}u, \sqrt{t_0}v)\in N(\ve,\O, V_1, V_2)$,
$$c_{\ve, \O, V_1, V_2}\leq E_{\ve, \Om, V_1, V_2} (\sqrt{t_0}u, \sqrt{t_0}v) \leq \sup_{t>0}
E_{\ve,\Om, V_1, V_2} [\sqrt{t}u, \sqrt{t}v]$$ and hence $c_{\ve,
\O, V_1, V_2}\leq m_\ve$. Therefore, we may complete the proof of
this Lemma.

\end{proof}

\section{Proofs of Theorem \ref{t1.1} and Theorem \ref{t1.2}}
\setcounter{equation}{0}

In this section, we prove Theorem \ref{t1.1} and Theorem
\ref{t1.2} by approximation argument. Fix a ball $\O=B_k$, where
$k$ is a large parameter tending to infinity. By
Theorem~\ref{t2.1}, each $c_{\ve, B_k, V_1, V_2}$ is attained by $
(u_k, v_k)$ a least energy solution of the following problem:
\begin{equation}\label{eq2.1}
\begin{cases}
\ve^2\triangle u(x)-V_1(x)u(x)+\mu_1 u^3 + \b u v^2=0\ \mathrm{in}
\  B_k,\\
\ve^2\triangle v(x)-V_2(x)v(x)+ \mu_2 v^3 + \b u^2 v=0\
\mathrm{in}
\  B_k,\\
u, v>0 \ \ \mathrm{in}
\  B_k, \ u=v=0\ \mathrm{on}\ \partial B_k.
\end{cases}
\end{equation}
By examining the argument in the proof of Theorem~\ref{t2.1}, we
may obtain the following estimates:
\begin{equation}\label{c8}
C_1\ve^N\leq\int_{B_k}u_{k}^4\leq C_2\ve^N,\quad
C_1\ve^N\leq\int_{B_k}v_{k}^4\leq C_2\ve^N\,,
\end{equation}
where $C_1$ and $C_2$ are positive constants independent of
$0<\ve\leq1$ and $k\geq1$. By the system~(\ref{eq2.1}) and
(\ref{c8}), we may derive that
\begin{equation}\label{2.1-1}
\int_{B_k} [ \ve^2 |\nabla u_{k}|^2+V_1 u_k^2 + \ve^2 |\nabla
v_k|^2 + V_2 v_k^2]\leq C_3\ve^N\,,
\end{equation} where $C_3$ is a positive constant independent of
$0<\ve\leq1$ and $k\geq1$. We may extend each $u_{k}$ and $v_k$
equal to 0 outside $B_k$, respectively. Then (\ref{2.1-1}) may
give \beq\label{he1}||u_{k}||_{H^1(\R^N)} +||v_{k}||_{H^1(\R^N)}
\leq C_4\ve^{N/2}\,,\eeq where $C_4$ is a positive constant
independent of $0<\ve\leq1$ and $k\geq1$.

Now we study the asymptotic behavior of $u_{k}, v_k$ as
$k\rightarrow\infty$. Due to (\ref{he1}), we obtain that as
$k\rightarrow\infty$, $u_{k}\rightharpoonup \bar{u},
v_{k}\rightharpoonup \bar{v}$, where $\bar{u}, \bar{v}\geq0$ and
$\bar{u}, \bar{v}\in H^1(\R^N)$. Moreover, the standard elliptic
regularity theorem may give that $(\bar{u}, \bar{v})$ is a
solution of the system \beq\label{system}
\left\{\begin{array}{llll} &\ve^2 \D \bar{u}-V_1 \bar{u} + \mu_1
\bar{u}^3 +
\b\bar{v}^2 \bar{u} &=0 &\mbox{in} \ \R^N\,, \\
&\ve^2 \D\bar{v}-V_2 \bar{v} + \mu_2 \bar{v}^3 + \b\bar{u}^2
\bar{v}&=0 \ &\mbox{in} \ \R^N\,.\end{array}\right. \eeq Then we
have the following lemma, whose proof is exactly same as those of
Theorem~3.3 in~\cite{lw3}.

\begin{lemma}\label{l2.3}
\hspace*{0.1cm}
\begin{itemize}
\item[(a)]As $k\rightarrow\infty$, $c_{\ve,B_k, V_1,
V_2}\rightarrow c_{\ve, \R^N, V_1, V_2}\,,$
 \item[(b)]If $\bar{u} \not \equiv 0, \bar{v} \not \equiv 0$, then
$(\bar{u},\bar{v})$ is a solution of (\ref{eq:1-4s}) and attains
$c_{\ve, \R^N, V_1, V_2}$, i.e. $(\bar{u},\bar{v})$ is a ground
state solution of (\ref{eq:1-4s}).
\end{itemize}
\end{lemma}
It remains to show that $ \bar{u} \not \equiv 0, \bar{v} \not
\equiv 0$. Note that if $\bar{u} \equiv 0$, then $\bar{v}$
satisfies
\begin{equation}
\ve^2 \Delta \bar{v}- V_2 \bar{v} + \mu_2 \bar{v}^3=0\,.
\end{equation}
Due to $\mu_2 \leq 0$, it is obvious that $\bar{v} \equiv 0$.
Therefore, we only need to exclude the case that $\bar{u}
\equiv\bar{v}\equiv 0$.

Suppose $ V(x) \equiv \l_1$ and $V_2 (x) \equiv \l_2$. Then by the
Maximum Principle and Moving Plane Method, both $u_k$ and $v_k$
are radially symmetric, strictly decreasing and satisfy
\begin{equation}{\label{bf2k}} \left\{\begin{array}{l}
\ve^2 \D u_k- \l_1  u_k+\mu_1 u_k^3 + \b u_k v_k^2 =0\ \ \ \mbox{in}\ \ B_k,\\
\ve^2 \D v_k- \l_2  v_k + \mu_2 v_k^3+ \b u_k^2 v_k =0\ \ \ \mbox{in}\ \ B_k \,,\\
u_k=u_k(r),v_k=v_k(r) >0\ \ \mbox{in}\ \ B_k,\\
u=v=0\ \ \mbox{on}\ \ \p B_k.
\end{array}\right.
\end{equation} Here we have used the fact that $\l_j>0$,
$\mu_j\leq 0\,, j=1, 2$ and $\beta>0$. Moreover, since the origin
$0$ is the maximum point of $u_k$ and $v_k$, then $\D u_k(0)\,, \D
v_k(0)\leq 0$ and $u_k(0)\,, v_k(0)>0$. Hence by (\ref{bf2k}), we
have
$$
 \b (v_k (0))^2 \geq -\mu_1 (u_k(0))^2 + \l_1, \ \ \b
(u_k (0))^2 \geq -\mu_2 (v_k(0))^2 + \l_2\,.
$$
Consequently, as $k \to +\infty$, \begin{eqnarray}\label{m3.2} \b
(v_0 (0))^2 &\geq& -\mu_1 (u_0(0))^2 + \l_1 \geq \l_1\,, \\
\nonumber \b (u_0 (0))^2 &\geq&  -\mu_2 (v_0(0))^2 + \l_2\geq
\l_2\,.\end{eqnarray} Here we have used the fact that $\mu_j\leq
0$ and $(u_k, v_k) \to (u_0, v_0)$ in $C_{loc}^2 (\R^N)$.
Therefore, (\ref{m3.2}) may imply that $ u_0 \not \equiv 0, v_0
\not \equiv 0$ and $(u_0, v_0) \in N(1, \R^N, \l_1, \l_2)$ is a
minimizer of $ c_{1, \R^N, \l_1, \l_2}$.

On the other hand, any minimizer of $c_{1, \R^N, \l_1, \l_2}$,
called $(U_0,V_0)$, must satisfy
\begin{equation}{\label{eq:3-20}}
\left\{ \begin{array}{l}
\D U_0 - \l_1 U_0 + \mu_1 U_0^3 + \b U_0 V_0^2=\ 0 \ \ \mbox{in} \ \ \R^N,\\
\D V_0 - \l_2 V_0 + \mu_2 V_0^3 + \b U_0^2 V_0 =\  0 \ \ \mbox{in}\ \ \R^N,\\
U_0, V_0 > 0, U_0, V_0 \in H^1 (\R^N). \end{array}\right.
\end{equation}
Due to $\b>0$, the problem~(\ref{eq:3-20}) is of cooperative
systems. By the moving plane method (cf.~\cite{T}), $(U_0, V_0)$
must be radially symmetric and strictly decreasing. This may
complete the proof of Theorem~\ref{t1.1}.

To finish the proof of Theorem \ref{t1.2}, we divide the proof into two cases as
follows:

\noindent {\bf Case 1:} either $ b_1^\infty = \infty$ or $
b_2^\infty =\infty$.

\begin{proof}

In this case, we note that
\begin{align}
c_{\ve, B_k, V_1, V_2}&=\frac{1}{4} \int_{B_k} \bigg[ \mu_1 u_{k}^4+2\beta u_{k}^2v_{k}^2+\mu_2v_{k}^4\bigg]\notag\\
&\leq C_3\ve^N\,,\notag\end{align} and \begin{align} c_{\ve, B_k,
V_1, V_2}&=\frac{1}{4} \int_{B_k} \bigg[ \ve^2 |\nabla u_{k}|^2+
V_1 u_k^2 + \ve^2 |\nabla v_k|^2 + V_2 v_k^2\bigg]
\notag\\
&\geq
C_4\ve^{N/2}\Bigg(\sqrt{\int_{B_k}u_{k}^4}+\sqrt{\int_{B_k}v_{k}^4}\Bigg)\,.\notag
\end{align}
Consequently,
\begin{align}\label{id6.2}
C_5\ve^N\leq c_{\ve, B_k, V_1, V_2}\leq C_6\ve^N,
\end{align}
where $C_5,C_6$ are independent of $\ve\leq1,\ k\geq1$. This
 gives
\begin{align}\notag
\int_{B_k} [ \ve^2 |\nabla u_{k}|^2+ V_1 u_k^2 + \ve^2 |\nabla
v_k|^2 + V_2 v_k^2] \leq C_7\ve^N.
\end{align}
By Sobolev's embedding (since $N \leq  3$),
\begin{equation}\label{id6.3}
\int_{B_k}u_{k}^6\leq C_8\ve^N,\quad \int_{B_k\cap\{|x|\geq
R\}}u_{k}^2\leq C_9\ve^N\cdot\frac{1}{\min\limits_{|x|\geq
R}V_1(x)}\,.
\end{equation}
Hence
\begin{align}\label{id6.4}
\int_{B_k\cap\{|x|\geq R\}}u_{k}^4&\leq\bigg(\int_{B_k\cap\{|x|\geq R\}}
u_{k}^2\bigg)^{1/2}\bigg(\int_{B_k\cap\{|x|\geq R\}}u_{k}^6\bigg)^{1/2
}\notag\\
&\leq C_{10} \ve^N \cdot\left(\frac{1}{\min\limits_{|x|\geq
R}V_1(x)}\right)^{1/2}.
\end{align}
By (\ref{c8}) and (\ref{id6.4}), we have
\begin{align}\label{id6.5}
\int_{B_k\cap\{|x|\leq
R\}}u_{k}^4\geq\Bigg(C_1-\frac{C_{10}}{\sqrt{\min\limits_{|x|\geq
R}V_1(x)}}\Bigg)\ve^N.
\end{align}
Thus if $u_{k}\rightharpoonup\overline{u}$, then $\overline{u}\geq0$
 and
\begin{align}\label{id6.6}
\int_{B_R
}\overline{u}^4\geq\Bigg(C_1-\frac{C_{10}}{\sqrt{\min\limits_{|x|\geq
R}V_1(x)}}\Bigg)\ve^N.
\end{align}
Due to $b_1^\infty=+\infty$, we may choose $R$ large enough such
that $C_1-\frac{C_{10}}{\sqrt{\min\limits_{|x|\geq
R}V_1(x)}}\geq\frac{1}{2}C_1$. Consequently, $\int_{B_R
}\overline{u}^4\geq\frac{1}{2}C_1\ve^N$ and hence
$\overline{u}\not\equiv 0$.\\
\end{proof}

\noindent
{\bf Case 2:} $ b_j^\infty<+\infty,\ j=1,2$

\begin{proof}

Suppose $\overline{u}\equiv\overline{v}\equiv0$. Then
\begin{equation}\label{id3.2}
u_{k}, v_k\rightarrow0\ \mathrm{in}\ \mathbb{C}_{loc}^2(\R^N).
\end{equation}
Let $M$ and $R$ be such that
\begin{equation}\label{id3.3}
|V_j(x)-b_j^\infty|<\frac{1}{M}\quad\hbox{ for }\: |x|\geq R\,.
\end{equation}
Let $\chi_R(x)$ be a smooth cut-off function such that
$\chi_R(x)=1$ for $|x|\leq R$, $\chi_R(x)=0$ for $|x|\geq 2R$. Now
we set
\begin{equation}\label{id3.4}
\widetilde{u}_{k}=u_{k}(1-\chi_R)\,,\quad\widetilde{v}_{k}=v_{k}(1-\chi_R)\,.
\end{equation}
Then we have
\begin{equation}\notag
\int_{\R^N}|\nabla \widetilde{u}_{k}|^2=\int_{\R^N}|\nabla
u_{k}|^2 -2\int_{\R^N}\nabla u_{k}\cdot\nabla(u_{k}\chi_R)
+\int_{\R^N}|\nabla (u_{k}\chi_R)|^2,
\end{equation} and
\begin{equation}\notag
\lim_{k\rightarrow+\infty}\Bigg(\Bigg|\int_{\R^N}\nabla
u_{k}\cdot\nabla(u_{k}\chi_R)\Bigg| +\int_{\R^N}|\nabla
u_{k}\chi_R|^2\Bigg)=0\,.
\end{equation}
Now we denote $o(1)$ as the terms that approach zero as
$k\rightarrow\infty$. Thus we can write
\begin{equation}\label{id3.5}
\int_{\R^N}|\nabla \widetilde{u}_{k}|^2=\int_{\R^N}|\nabla u_{k}|^2+o(1).
\end{equation}
Similarly,
\begin{equation}\notag
\int_{\R^N}|\nabla \widetilde{v}_{k}|^2=\int_{\R^N}|\nabla v_{k}|^2+o(1), \int_{\R^N}V_1\widetilde{u}_{k}^p=\int_{\R^N}V_1 u_{k}^p+o(1), \int_{\R^N}V_2\widetilde{v}_{k}^p=\int_{\R^N}V_2 v_{k}^p+o(1)
\end{equation}
for all $ 2\leq p \leq 6$. Hence $E_{\ve,B_k, V_1,
V_2}[u_{k},v_{k}]=c_{\ve, B_k, V_1, V_2}=E_{\ve, B_k, V_1,
V_2}[\widetilde{u}_{k},\widetilde{v}_{k}]+o(1)$. Moreover,
\begin{align}\label{id3.6}
&\int_{\R^N} [\ve^2 |\nabla \widetilde{u}_{k}|^2+b_1^\infty
\widetilde{u}_{k}^2 +\ve^2 |\nabla \widetilde{v}_{k}|^2+
b_2^\infty \widetilde{v}_{k}^2 ]
\\
 - & \int_{\R^N} [  \mu_1\widetilde{u}_{k}^4+2\beta \widetilde{u}_{k}^2\widetilde{u}_{k}^2 +\mu_1\widetilde{v}_{k}^4]
\notag\\
=&\int_{\R^N}(b_1^\infty-V_1(x))\widetilde{u}_{k}^2+\int_{\R^N}(b_2^\infty-V_2(x))\widetilde{v}_{k}^2 +o(1)\notag\\
=&O\Big(\frac{1}{M}\int_{\R^N}(\widetilde{u}_{k}^2 + \widetilde{v}_k^2) \Big)+o(1)\notag\\
=&O\Big(\frac{1}{M}\Big)+o(1),\ j=1,2\,.\notag
\end{align}
Similarly, we have
\begin{equation}\label{id3.6-1}
\int_{\R^N} [ 2 \b \widetilde{u}_k^2 \widetilde{v}_k^2 +\mu_1
\widetilde{u}_k^4+ \mu_2 \widetilde{v}_k^4] = \int_{\R^N} [ 2 \b
u_k^2 v_k^2 +\mu_1 u_k^4+ \mu_2 v_k^4] + o(1) C \ve^N.
\end{equation}
Hence by (\ref{id3.6}), (\ref{id3.6-1}) and (\ref{eq:2-25}) of
Claim~1 in Theorem~\ref{t2.1}, we see that the unique critical
point $\widetilde{t}$ of the function $E_{\ve, \R^N, b_1^\infty,
b_2^\infty}[\sqrt{t}\widetilde{u}_{k},\sqrt{t}\widetilde{v}_{k}]$
satisfies
\begin{equation}\label{id3.7}
|\widetilde{t}-1|=O\Big(\frac{1}{M}\Big)+o(1)\,,
\end{equation}
which yields
\begin{align}
E_{\ve, \R^N, b_1^\infty, b_2^\infty}
\Big[\sqrt{\widetilde{t}}\widetilde{u}_{k},\sqrt{\widetilde{t}}\widetilde{v}_{k}\Big]
=&E_{\ve, \R^N, b_1^\infty, b_2^\infty}[\widetilde{u}_{k},\widetilde{v}_{k}]+O\Big(\frac{1}{M}\Big)+o(1)\notag\\
=&E_{\ve, \R^N, V_1, V_2}[\widetilde{u}_{k},\widetilde{v}_{k}]+O\Big(\frac{1}{M}\Big)+o(1)\notag\\
=&E_{\ve, \R^N, V_1, V_2}[u_{k},v_{k}]+O\Big(\frac{1}{M}\Big)+o(1)\notag\\
=&c_{\ve, B_k, V_1, V_2}+O\Big(\frac{1}{M}\Big)+o(1).\notag
\end{align}
On the other hand,
\begin{equation}\label{id3.8}
\Big(\sqrt{\widetilde{t}}\widetilde{u}_{k},\sqrt{\widetilde{t}}\widetilde{v}_{k}\Big)\in
N (\ve, \R^N, b_1^\infty, b_2^\infty)
\end{equation}
and then
\begin{equation}\label{id3.9}
E_{\ve, \R^n, b_1^\infty, b_2^\infty}
\Big[\sqrt{\widetilde{t}}\widetilde{u}_{k},\sqrt{\widetilde{t}}\widetilde{v}_{k}\Big]
\geq c_{\ve, \R^N, b_1^\infty, b_2^\infty}
\end{equation} Consequently, $c_{\ve, \R^N, b_1^\infty, b_2^\infty}\leq
c_{\ve, B_k, V_1, V_2}+O\big(\frac{1}{M}\big)+o(1)$. Letting
$M\rightarrow+\infty$ and $k\rightarrow+\infty$, we obtain
$c_{\ve, \R^N, b_1^\infty, b_2^\infty}\leq c_{\ve, \R^N, V_1,
V_2}$ which may contradict with (\ref{c7}). Therefore, we may
complete the proof of Theorem~\ref{t1.2}.
\end{proof}

\section{Proof of Theorem~\ref{t1.3}}
\setcounter{equation}{0}

\ \ \ \ In this section, we study the asymptotic behavior of
$(u_{\ve},v_{\ve})$ as $\ve\rightarrow 0$. Firstly, the energy
upper bound is stated as follows:
\begin{lemma}\label{lem5.1}
For $\beta>0$ and $0<\ve<<1$,
\begin{equation}\label{id5.1}
c_{\ve, \R^N, V_1, V_2} \leq\ve^N [\inf\limits_{x\in\R^N} c_{1,
\R^N, V_1(x),V_2(x) } +o(1)].
\end{equation}
\end{lemma}
\begin{proof}[\bf P{\small ROOF}]
Fix a point $x_0 \in \R^N$. Let $(U_0, V_0)$ be a minimizer of $
c_{1, \R^N, V_1 (x_0), V_2 (x_0)}$. We set $ u(x)= U_0
(\frac{x-x_0}{\ve}), v(x)= V_0 (\frac{x-x_0}{\ve})$ and then use
(\ref{eq:2-25s}) to compute the upper bound of $ c_{\ve, \R^N,
V_1, V_2}$. Due to $ c_{\ve, \R^N, \l_1, \l_2} = \ve^N c_{1, \R^N,
\l_1, \l_2}$, the  rest of the proof is simple and thus omitted.
\end{proof}

Let $u_{\ve}(P^\ve)=\sup\limits_{x\in\R^N}u_{\ve}(x)$ and
$v_{\ve}(Q^\ve)=\sup\limits_{x\in\R^N}v_{\ve}(x)$. We want to
claim that $\sup\limits_{\ve>0}(|P^\ve|+|Q^\ve|)<+\infty$. To this
end, we need to show that both $u_{\ve}$ and $v_{\ve}$ are
uniformly bounded. In fact, as for the proof of (\ref{id6.3}), we
have
\begin{equation}\label{id6.7}
\int_{\R^N}(u_{\ve}^q + v_\ve^q)\leq c\ve^N,\ 2\leq q\leq6.
\end{equation}
The equation of $u_{\ve}$  gives
\begin{align}
\ve^2\triangle u_{\ve}=&V_1u_{\ve}-\mu_1u_{\ve}^3
-\beta u_{\ve}v_{\ve}^2\notag\\
\geq&-\beta v_{\ve}^2u_{\ve}\notag\\
=&-C(x)u_{\ve}\quad\hbox{ in }\:\mathbb{R}^N\,.\notag
\end{align}
Let $\widetilde{U}_{\ve}(y)=u_{\ve}(\ve\,y)$, and $C_\ve(y)=
C(\ve\,y)$. Then
\begin{equation}\label{id6.8}
\triangle\widetilde{U}_{\ve}+C_\ve(y)\widetilde{U}_{\ve} \geq
0\quad\hbox{ in }\:\mathbb{R}^N\,,\quad\hbox{ and }\: C_\ve\in
L^3(\R^N)\,.
\end{equation}
By the subsolution estimate (Theorem~8.17 of~\cite{GT})
\begin{equation}\label{id6.9}
|\widetilde{U}_{\ve}(y)|\leq
C\Bigg(\int_{B(y,1)}|\widetilde{U}_{\ve}|^2\Bigg)^{1/2},
\end{equation}
where $C>0$ is independent of $\ve$. Hence by (\ref{id6.7}) and
(\ref{id6.9}), we see that $||\widetilde{U}_{\ve}||_{L^\infty}\leq
C$ and hence $0<u_{\ve}\leq C$. Similarly, we may obtain
$0<v_{\ve}\leq C$.

\noindent {\bf Claim~3:}~~{\it If $ |P^\ve| \to +\infty$, then $
b_1^\infty <+\infty$.} Suppose $ b_1^\infty=+\infty$. Since
$P^\ve$ is a local maximum point of $u_{\ve}$, then $\triangle
u_{\ve}(P^\ve)\leq 0$. Hence by the equation of $u_{\ve}$, we may
obtain
\begin{equation}
V_1(P^\ve)u_{\ve}(P^\ve)-\mu_1 u_{\ve}^3(P^\ve)-\beta u_{\ve
}(P^\ve)v_{\ve}^2(P^\ve)=\ve^2\triangle u_{\ve}(P^\ve)\leq0,\notag
\end{equation}
which implies that
\begin{equation}\label{id6.10}
V_1(P^\ve)\leq\beta v_\ve^2 (P^\ve)\leq C,
\end{equation}
and hence
\begin{equation}\label{id6.11}
|P^\ve|\leq C_0\,.\end{equation} Therefore, we may complete the
proof of Claim~3. Moreover, we may also claim that $b_2^\infty
<+\infty$. In fact, suppose $b_2^\infty=+\infty$. Set $ U_\ve(y):=
u_\ve (P^\ve+\ve y), V_\ve (y):= v_\ve (P^\ve +\ve y) $. Then
$U_\ve \to U_0$ in $C^2_{loc} (\R^N)$ and $ V_\ve \to V_0$  in
$C^2_{loc} (\R^N)$, where $(U_0, V_0)$ satisfies
\begin{equation}
\label{uv00}
\D  U_0 - b_1^\infty U_0 + \mu_1 U_0^3 +\b U_0 V_0^2=0 \ \mbox{in } \ \R^N.
\end{equation}
Hence by (\ref{id6.10}), we may obtain $V_0(0)>0$, and then  $ V_0
\not \equiv 0$. This implies that
\begin{eqnarray*}
c_{\ve, \R^N, V_1, V_2} & = & \frac{1}{4}\int_{\R^N} [ \ve^2
|\nabla u_\ve|^2 +V_1 u_\ve^2 +\ve^2 |\nabla v_\ve|^2 + V_2
v_\ve^2]
 \\
 &  \geq  & \frac{1}{4}\int_{|x|>R}
 [ \ve^2 |\nabla u_\ve|^2 +V_1 u_\ve^2 +\ve^2 |\nabla v_\ve|^2 + V_2 v_\ve^2]
\\
& \geq &  \frac{1}{4} \int_{|x|>R} V_2 v_\ve^2 \\
& \geq & C \ve^N \left[ \inf_{|x|>R} V_2 (x)\right]
\end{eqnarray*}
which contradicts with (\ref{id5.1}). Here we have used the
hypothesis that $b_2^\infty=+\infty$. Thus we may assume that
$b_1^\infty <+\infty$ and $b_2^\infty <\infty$. As before,
$(U_\ve, V_\ve)$ converges to $(U_0, V_0)$ satisfying
\begin{equation}
\label{uv0} \D  U_0 - b_1^\infty U_0 + \mu_1 U_0^3 +\b U_0
V_0^2=0, \ \D  V_0 - b_2^\infty V_0 + \mu_1 V_0^3 +\b V_0 U_0^2=0
\ \mbox{in } \ \R^N\,.
\end{equation}
Then again $V_0 \not \equiv 0$ since otherwise, $ (U_0,V_0)\equiv
(0,0)$ which is impossible. Moreover,
\begin{eqnarray*}
c_{\ve, \R^N, V_1, V_2} & = & \frac{1}{4}\int_{\R^N} [ \ve^2
|\nabla u_\ve|^2 +V_1 u_\ve^2 +\ve^2 |\nabla v_\ve|^2 + V_2
v_\ve^2]
 \\
 &  \geq  & \frac{1}{4}\int_{|x|>R} [ \ve^2 |\nabla u_\ve|^2 +V_1 u_\ve^2 +\ve^2 |\nabla v_\ve|^2 + V_2 v_\ve^2]
\\
& \geq & \ve^N   \frac{1}{4}\int_{\R^N} [ |\nabla U_0|^2 +b_1^\infty U_0^2 +|\nabla V_0|^2 + b_2^\infty V_0^2] + o(\ve^N) \\
& \geq & \ve^N [ c_{1, \R^N, b_1^\infty, b_2^\infty} + o(1)]
\end{eqnarray*}
which may contradict with (\ref{id5.1}). Therefore, we complete
the proof of $\ds\sup_{\ve>0}|P^\ve |+|Q^\ve| <+\infty$.

Let $(P^\ve, Q^\ve)\rightarrow(P^0, Q^0)$. As before, $(U_\ve,
V_\ve)= (u_{\ve} (P^\ve +\ve y), v_\ve (P^\ve +\ve y)) \to (U_0,
V_0)$, where $(U_0, V_0)$ satisfies
\begin{equation}\notag
\begin{cases}
\triangle U-V_1(P^0)U+\mu_1 U^3+\beta U V^2=0\quad\hbox{ in }\:\R^N\,,\\
\triangle V-V_2(P^0)V+\mu_2 V^3+\beta U^2 V=0\quad\hbox{ in
}\:\R^N\,.
\end{cases}
\end{equation}
Then by the strong Maximum Principle, $U_0, V_0>0$. Furthermore,
we have
\begin{equation}\notag
\lim_{\ve\rightarrow0}\ve^{-N}c_{\ve, \R^N, V_1, V_2}\geq c_{1,
\R^N, V_1(P^0),V_2( P^0)}.
\end{equation}
Hence by Lemma~\ref{lem5.1},
\begin{equation}\notag
c_{1, \R^N, V_1(P^0),V_2(P^0)} \leq\inf_{x\in\R^N} c_{1, \R^N,
V_1(x),V_2(x)}\,,
\end{equation} i.e.
$c_{1, \R^N, V_1(P^0),V_2(P^0)}=\ds\inf_{x\in\R^N} c_{1, \R^N,
V_1(x),V_2(x)}\,.$

It remains to show that $\frac{|P^\ve-Q^\ve|}{\ve}\rightarrow 0$.
In fact, if $\frac{|P^\ve-Q^\ve|}{\ve}\to +\infty$, then similar
arguments may give
\begin{equation}\notag
\lim_{\ve\rightarrow0}\ve^{-N}c_{\ve, \R^N, V_1, V_2}\geq c_{1,
\R^N, V_1(P^0),V_2( P^0)}+ c_{1, \R^N, V_1(Q^0),V_2( Q^0)} \geq 2
\inf_{x \in \R^N} c_{1, \R^N, V_1 (x), V_2 (x)}
\end{equation}
which is impossible. On the other hand, if
$\frac{|P^\ve-Q^\ve|}{\ve}\to c \not = 0$, then $U_0$ and $V_0$
may have different maximum points. This may contradict with the
fact that both $U_0$ and $V_0$ are radially symmetric and strictly
decreasing. Thus $\frac{|P^\ve-Q^\ve|}{\ve}\rightarrow 0$. The
uniqueness of $P^\ve,Q^\ve$ may follow from Claim~8 of~\cite{lw1}.
Therefore, we may complete the proof of Theorem~\ref{t1.3}.

\section{ Proof of Theorem \ref{t1.4}}
\setcounter{equation}{0}

In this section, we follow the same ideas of~\cite{lw1} to prove
Theorem~\ref{t1.4}. As for the proof of Lemma~4.2 in~\cite{lw1},
the upper bound of $c_{\ve, \O, \l_1, \l_2}$ is given by
\begin{lemma}
For $\b > \sqrt{\mu_1 \mu_2}$,
\begin{equation}
{\label{eq:4-17}} c_{\ve, \O, \l_1, \l_2}  \leq   \ve^N \Biggl\{
c_{1, \R^N, \l_1, \l_2}+  c_1\,e^{-2\sqrt{\l_1}(1-\s) R_\ve} +
c_2\,e^{-2\sqrt{\l_2}(1-\s)R_\ve } \Biggl\}\,,
\end{equation}
where $R_\ve = \frac{1}{\ve} \ds\max_{P \in \Omega} d(P,
\partial \Omega)$ and $c_j$'s are positive constants. \label{upper}
\end{lemma}
\noindent Furthermore, the asymptotic behavior of $(u_\ve,
v_\ve)$'s can be summarized as follows:
\begin{lemma}
For $\ve$ sufficiently small, $u_\ve$ has only one local maximum
point $P_\ve$ and $v_\ve$ has only one local maximum point $Q_\ve$
such that
\begin{equation}{\label{eq:5-5}}
\frac{d(P_\ve,\p\O)}{\ve} \to +\i,\quad \frac{d(Q_\ve,\p\O)}{\ve}
\to +\i,\quad \frac{|P_\ve-Q_\ve|}{\ve}\to 0.
\end{equation}
Let $U_\ve(y) :=\  u_\ve (P_\ve +\ve y)$, $V_\ve(y) :=\ (Q_\ve+\ve
y)$. Then $(U_\ve, V_\ve)\to (U_0,V_0)$, where $(U_0,V_0)$ is a
least-energy solution of (\ref{eq:1-4g}). Moreover,
\begin{equation}
\label{uvdecay} \ve \left|\bd u_\ve\right|+|u_\ve|\leq C
 e^{-\sqrt{\l_1}(1-\s)\frac{|x-P_\ve|}{\ve}},
 \quad \ve \left|\bd v_\ve\right|+|v_\ve|\leq C
 e^{-\sqrt{\l_2}(1-\s)\frac{|x-Q_\ve|}{\ve}}\,.
\end{equation}
\end{lemma}

Now we want to complete the proof of Theorem~\ref{t1.4}. We may
assume that, passing to a subsequence, that $P_\ve$ (or $Q_\ve)
\to x_0\in \bar \O$. Thus
$$d_\ve =\  d(P_\ve,\p\O) \to d_0 :=\  d(x_0,\p\O),\ \ \mbox{as}\ \ \ve\to
 0.$$
Note that $d_0$ may be zero. Given $\s>0$ a small constant, we may
choose $d'_0>0$ and $\s'>0$ slightly smaller than $\s$ such that
$$\mbox{vol}(B(x_0,d'_0))=\  \mbox{vol} (\O\cap B(x_0,d_0+\s))
\quad\hbox{ and }\quad d'_0<d_0+\s'\,.$$ Besides, we may set
$\eta_\ve$ as a $C^\i$ cut-off function such that
$$\left\{\begin{array}{lll} &\eta_\ve(s) =\ 1\ \ &\mbox{for}\ \ 0\leq s\leq d_\ve+\s'\,,\\
&\eta_\ve(s) =\ 0 \ \ &\mbox{for}\ \ s>d_\ve+\s\,, \\
&0\leq\eta_\ve\leq 1\,, &|\eta'_\ve|\leq C\,.\end{array}\right.$$
Let $\tilde u_\ve(x)=\ u_\ve \eta_\ve(|P_\ve-x|)$ and $\tilde
v_\ve(x) =\ v_\ve\eta_\ve(|Q_\ve-x|)$. Then we have
\begin{equation}
\lim_{\ve \to 0} \ve^{-N}  \int_\O [ 2 \b \tilde{u}_\ve^2
\tilde{v}_\ve^2 + \mu_1 \tilde{u}_\ve^4 +\mu_2 \tilde{v}_\ve^4 ] =
\int_{\R^N} [ 2 \b U_0^2 V_0^2 +\mu_1 U_0^4 + \mu_2 V_0^4]>0\,.
\end{equation}
Hence
\[  \int_\O [ 2 \b \tilde{u}_\ve^2 \tilde{v}_\ve^2 + \mu_1 \tilde{u}_\ve^4
+\mu_2 \tilde{v}_\ve^4 ]>0\,,\] as $\ve$ sufficiently small.

By the decay estimate~(\ref{uvdecay}) and Lemma~\ref{cnew},  we
obtain that \begin{eqnarray}\label{6.4-1} c_{\ve, \O, \l_1, \l_2}
&\geq& E_{\ve, \Omega, \l_1, \l_2} [tu_\ve,t v_\ve]  \\
&\geq& E_{\ve, \tilde{\Omega}, \l_1, \l_2} [t\tilde u_\ve, t\tilde
v_\ve] -\ve^N \exp \left[-\frac{2\sqrt{\l_1}}{ \ve}
(d_\ve+\s')\right] -\ve^N \exp \left[-\frac{2\sqrt{\l_2}}{ \ve}
(d_\ve+\s')\right]\nonumber
\end{eqnarray} for all $t\in [0,2]$, where
$\tilde{\Omega}= \Omega \cap B(x_\ve, d_\ve +\s)$ and $x_\ve$ can
be $P_\ve$ or $Q_\ve$. Let $R_\ve=\ \frac{d'_\ve}{\ve}$, where
$d'_\ve$ is chosen such that
$$\mbox{vol} (B(0,d'_\ve)) =\  \mbox{vol} (\O\cap B(x_\ve,d_\ve+\s)).$$
Using Schwartz's symmetrization, we have
$$\int_{B(0, d'_\ve)} (\tilde u_\ve^*)^2 (\tilde v_\ve^*)^2 \geq \int_{\tilde{\O}} \tilde u_\ve^2
\tilde v_\ve^2$$ and then
\begin{equation}
\label{uvin} \int_{B(0, d'_\ve)} [ 2 \b (\tilde u_\ve^*)^2 (\tilde
v_\ve^*)^2 + \mu_1 (\tilde{u}_\ve^*)^4 + \mu_2
(\tilde{v}_\ve^*)^4] \geq \int_{\tilde{\O}} [ 2 \b \tilde u_\ve^2
\tilde{v}_\ve^2 +\mu_1 \tilde{u}-\ve^4 +\mu_2 \tilde{v}_\ve^4]
>0\,.
\end{equation}
Thus \beq\label{6.6} E_{\ve, B(0, d'_\ve), \l_1, \l_2} [t\tilde
u_\ve^*, t \tilde v_\ve^*] \leq E_{\ve, \tilde{\Om}, \l_1, \l_2}
[t \tilde u_\ve, t \tilde
 v_\ve]\,,\quad\forall\, t\in [0,2]\,.\eeq Here we have used the fact that $\beta >0$.

By~(\ref{uvin}) and Claim~1 of Theorem~\ref{t2.1}, there exists
$t^*\in (0,2]$ such that
$$E_{\ve, B(0, d'_\ve), \l_1, \l_2} [t^* \tilde u_\ve^*, t^* \tilde v_\ve^*]
\geq E_{ \ve, B(0, d'_\ve), \l_1, \l_2}[t \tilde u_\ve^*, t\tilde
v_\ve^*]\,, \quad\forall t \geq 0\,.$$ Then by (\ref{6.4-1}) and
(\ref{6.6}),
\begin{eqnarray*}
& \ & E_{\ve, B(0, d'_\ve), \l_1, \l_2} [t^* \tilde u_\ve^*, t^* \tilde v_\ve^*] \\
& \leq &  E_{\ve, \tilde{\Om}, \l_1, \l_2} [t^* \tilde u_\ve, t^* \tilde v_\ve]\\
& \leq & c_{\ve, \O, \l_1, \l_2} + \ve^N \exp
\left[-\frac{2\sqrt{\l_1}}{\ve} (d_\ve+\s')\right] + \ve^N \exp
\left[-\frac{2\sqrt{\l_2}}{\ve} (d_\ve+\s')\right],
\\
 & \ & E_{\ve, B(0, d'_\ve), \l_1, \l_2}[t^* \tilde u_\ve^*, t^* \tilde v_\ve^*] \\
& =\  & \sup_{t>0} E_{\ve, B(0, d'_\ve), \l_1, \l_2} [t\tilde u_\ve^*, t\tilde v_\ve^*]\\
& \geq & \ve^N  \inf_{u, v \geq 0, \atop{u \not \equiv 0, v \not
\equiv 0,
\atop{(u, v)\in N(1,R_\ve, \l_1, \l_2) } } } E_{1, B_{R_\ve}, \l_1, \l_2} [u,v]\\
& \geq & \ve^N \Biggl\{ c_{1, \R^N, \l_1, \l_2} +  c_3\,\exp
\Biggl[ -\frac{2 (1+\s)\sqrt{\l_1}}{ \ve}
(d_\ve+o(1))\Biggl]\Biggl\} \\ & & + \ve^N\, c_4\,\exp \Biggl[
-\frac{2(1+\s)\sqrt{\l_2}}{ \ve} (d_\ve+o(1))\Biggl]\,,
\end{eqnarray*}
where $c_j$'s are positive constants. Here the last inequality may
follow from Lemma~\ref{upper} and Theorem~4.1 of \cite{lw1}. Thus
\beq \label{lowerbound2} c_{\ve, \O, \l_1, \l_2} \geq \ve^N
\left\{\begin{array}{lll} & c_{1, \R^N, \l_1, \l_2} +  c_3\,\exp
\Biggl[ -\frac{2 (1+\s)\sqrt{\l_1}}{ \ve} (d_\ve+o(1))\Biggl] \\
& + c_4\,\exp \Biggl[ -\frac{2(1+\s)\sqrt{\l_2}}{ \ve}
(d_\ve+o(1))\Biggl]\end{array} \right\}\,. \eeq
Combining the
lower and upper bound of $c_{\ve, \O, \l_1, \l_2}$, we obtain
\[  c_3\,\exp \Biggl[ -\frac{2 (1+\s)\sqrt{\l_1}}{ \ve} (d_\ve+o(1))\Biggl]
+ c_4\,\exp \Biggl[ -\frac{2(1+\s)\sqrt{\l_2}}{ \ve}
(d_\ve+o(1))\Biggl]
 \]
\[
\leq  c_1\,\exp \Biggl[ -\frac{2 (1-\s) \sqrt{\l_1}}{ \ve}
(d_0+o(1))\Biggl] +  c_2\,\exp \Biggl[ -\frac{2
(1-\s)\sqrt{\l_2}}{ \ve} (d_0+o(1))\Biggl]\,.
 \]
This then shows that $d(P_\ve,,\p\O), d(Q_\ve,\p\O) \to \ds\max_{P
\in\O} d(P,\p\O)$ since $  |P_\ve -Q_\ve| \to 0$. \qed

\end{document}